\documentclass[a4paper,12pt]{amsart}
\usepackage{amssymb}
\usepackage{graphicx}
\usepackage{hyperref}
\usepackage{amsmath}
\usepackage{amsthm}
\usepackage{amssymb}
 
\hypersetup{colorlinks =true,citecolor=blue,linkcolor=blue,filecolor=blue,urlcolor=blue,citecolor=blue}
\usepackage{ifthen}
\nonstopmode \numberwithin{equation}{section}
\newtheorem{thm}{Theorem}
\newtheorem{cor}{Corollary}
\newtheorem{lem}{Lemma}

\newtheorem{conj}{Conjecture}

\theoremstyle{definition}

\newtheorem{prob}[equation]{Problem}

\newenvironment{rem}{%
	\bigskip
	\noindent \textsl{{\sl Remark. }}}{\bigskip}
\newenvironment{rems}{%
	\bigskip
	\noindent \textsl{{\sl Remarks. }}}{\bigskip}

\newcounter {own}
\def\theown {\thesection       .\arabic{own}}

\newenvironment{pf}[1][]{%
	\vskip 3mm
	\noindent
	\ifthenelse{\equal{#1}{}}%
	{{\slshape Proof. }}%
	{{\slshape #1.} }%
}%
{\qed\bigskip}

\newcounter{alphabet}

\newcommand{\IR}{{\mathbb R}}
\newcommand{\ID}{{\mathbb D}}

\newcommand{\IC}{{\mathbb C}}

\newcommand{\D}{{\mathbb D}}

\newcommand{\dist}{{\operatorname{dist}}}

\def\be{\begin{equation}}
	\def\ee{\end{equation}}

\newcommand{\bee}{\begin{enumerate}}
	\newcommand{\eee}{\end{enumerate}}

\newcommand{\blem}{\begin{lem}}
	\newcommand{\elem}{\end{lem}}
\newcommand{\bthm}{\begin{thm}}
	\newcommand{\ethm}{\end{thm}}
\newcommand{\bcor}{\begin{cor}}
	\newcommand{\ecor}{\end{cor}}
\newcommand{\beg}{\begin{examp}}
	\newcommand{\eeg}{\end{examp}}
\newcommand{\begs}{\begin{examples}}
	\newcommand{\eegs}{\end{examples}}
\newcommand{\bdefe}{\begin{defin}}
	\newcommand{\edefe}{\end{defin}}
\newcommand{\bprob}{\begin{prob}}
	\newcommand{\eprob}{\end{prob}}
\newcommand{\bei}{\begin{itemize}}
	\newcommand{\eei}{\end{itemize}}

\newcommand{\bcon}{\begin{conj}}
	\newcommand{\econ}{\end{conj}}
\newcommand{\bcons}{\begin{conjs}}
	\newcommand{\econs}{\end{conjs}}
\newcommand{\bprop}{\begin{propo}}
	\newcommand{\eprop}{\end{propo}}
\newcommand{\br}{\begin{rem}}
	\newcommand{\er}{\end{rem}}
\newcommand{\brs}{\begin{rems}}
	\newcommand{\ers}{\end{rems}}
\newcommand{\bo}{\begin{obser}}
	\newcommand{\eo}{\end{obser}}
\newcommand{\bos}{\begin{obsers}}
	\newcommand{\eos}{\end{obsers}}
\newcommand{\bpf}{\begin{pf}}
	\newcommand{\epf}{\end{pf}}
\newcommand{\ba}{\begin{array}}
	\newcommand{\ea}{\end{array}}
\newcommand{\beq}{\begin{eqnarray}}
	\newcommand{\beqq}{\begin{eqnarray*}}
		\newcommand{\eeq}{\end{eqnarray}}
	\newcommand{\eeqq}{\end{eqnarray*}}

\newcounter{minutes}
\divide\time by 60
\newcounter{hours}
\multiply\time by 60 \addtocounter{minutes}{-\time}
\begin{document}
	\title{Length Distortion Of Curves Under Meromorphic Univalent Mappings}
	\begin{center}
		{\tiny \texttt{FILE:~\jobname .tex,
				printed: \number\year-\number\month-\number\day,
				\thehours.\ifnum\theminutes<10{0}\fi\theminutes}
		}
	\end{center}
	\author{Bappaditya Bhowmik${}^{~\mathbf{*}}$}
	\address{Bappaditya Bhowmik, Department of Mathematics,
		Indian Institute of Technology Kharagpur, Kharagpur - 721302,
		India.} \email{bappaditya@maths.iitkgp.ac.in}
	\author{Deblina Maity}
	\address{Deblina Maity, Department of Mathematics,
		Indian Institute of Technology Kharagpur, Kharagpur - 721302,
		India.}  \email{deblinamaity1997@gmail.com}
	
\begin{abstract}
	Let $f$ be a conformal (analytic and univalent) map defined on the open unit disk $\D$ of the complex plane $\IC$ that is continuous on the semi-circle $\partial \D^{+}=\{z\in\IC:|z|=1, {\rm{Im}}\,z>0\}$. The existence of a uniform upper bound for the ratio of the length of the image of the horizontal diameter $(-1,1)$ to the length of the image of $\partial \D^{+}$ under $f$ was proved by Gehring and Hayman. In this article, at first, we generalize this result by introducing a simple pole for $f$ in $\D$ and considering the ratio of the length of the image of the vertical diameter $I=\{z: {\rm{Re}}\,z=0; ~|{\rm{Im}}\,z|<1\}$ to the length of the image of the semi-circle $C'=\{z: |z|=1;~ {\rm{Re}}\,z<0\}$ under such $f$. Finally, we further generalize this result by replacing the vertical diameter $I$ with a hyperbolic geodesic symmetric with respect to the real line, and by replacing $C'$ with the corresponding arc of the unit circle passing through the point $-1$.
\end{abstract}
\thanks{The first author would like to thank SERB, India for its financial support through Core Research Grant
	(Ref. No.- CRG/2022/001835)}
\thanks{The second author of this article would like to thank
	the Prime Minister Research Fellowship of the Government of India (Grant ID: 2403440) for its financial support.}

\maketitle
\pagestyle{myheadings} \markboth{B. Bhowmik, D. Maity}{Length Distortion Of Curves Under Meromorphic Univalent Mappings}

\section{Introduction}
\noindent We will use the following notations throughout the article:
\vspace{0.2cm}
\\
$\IC:\text{The finite complex plane} $\\
$\hat{\IC}:\IC\cup\{\infty\}$\\
$\D:\{z\in \IC : |z|<1\}$\\
$\partial \D:\{z\in\IC : |z|=1\}$\\
$\partial \D^{+}:\{z\in\partial \D:~{\rm{Im}}\,z>0\}$\\
$ \mathbb{H}:\{z\in \IC : {\rm{Im}}\,z>0\}$\\
$ {\IR}^{+}:\{z\in \IC : {\rm{Im}}\,z=0;~ {\rm{Re}}\,z>0\}$\\
$  \mathbb{I}^{+}:\{z\in\mathbb{H} : {\rm{Re}}\,z=0\}$\\
$ I:\{z\in\D: {\rm{Re}}\,z=0; ~|{\rm{Im}}\,z|<1\}$\\
$C':\{ z\in \partial \D:~ {\rm{Re}}\,z<0\}$.
\vspace{0.2cm}

In this article, by a conformal mapping, we mean an analytic and univalent mapping. Let $f$ be a conformal mapping of $\ID$ that remains continuous on $\partial \D$. In \cite{FR}, Fej\'er and Riesz proved that, if $f$ maps $\partial \D$ onto a curve of length $L$, then the diameter $(-1,1)$ will be mapped onto a curve whose length does not exceed ${L}/{2}$. Later, Piranian conjectured (see \cite{GH}) the semi-circle version of the Fej\'er-Riesz theorem as stated below.
\vspace{0.2cm}

\noindent\textbf{Piranian's Conjecture.}
\textit{Suppose $f$ is a conformal mapping of $\mathbb{D}$ that remains continuous on the semi-circle $\partial \D^{+}$ and maps this semi-circle onto a curve of length $L$. Then $f$ maps the diameter $(-1,1)$ onto a curve whose length does not exceed $AL$, where $A$ is an absolute constant.}
\vspace{0.2cm}

The Fej\'er-Riesz theorem is a consequence of the following inequality which is valid for any function $f$ analytic in $\D\cup\partial\D$:
\be\nonumber
\int_{-1}^{1} |f(x)|dx \leq \frac{1}{2} \int_{-\pi}^{\pi} |f(e^{i\theta})|d\theta.
\ee
Unfortunately, this inequality is not true when considering its semi-circle version. However, Piranian's conjecture was ultimately validated. In \cite{GH}, Gehring and Hayman replaced  $\mathbb{D}$, $\partial {\D}^{+}$ and the horizontal diameter $(-1,1)$ with  $\mathbb{H}$, ${\IR}^{+}$ and $\mathbb{I}^{+}$, respectively, and proved the upper half-plane version of the Piranian's conjecture, which we state below (see \cite[Theorem~1]{GH}):
\vspace{0.2cm}

\noindent\textbf{Theorem A.}
\textit{If $f$ is a conformal mapping of $\mathbb{H}$ which remains continuous on the positive x-axis and maps this half-axis onto a curve of length $L$, then $f$ maps the positive y-axis onto a curve of length $L'$, where 
	$$
	L'\leq AL.
	$$
	Here, $A$ is an absolute constant with $\pi \leq A < 74$.}
\vspace{0.2cm}

In order to present the generalized Gehring-Hayman theorem, we first introduce the notion of the hyperbolic metric in $\D$. The hyperbolic metric in the unit disk $\D$ is defined by 
$$
\lambda_{\D}(z_1,z_2)=\underset{C}{\min}\int_{C}\frac{|dz|}{1-|z|^2}~~\mbox{for}~~z_1,z_2\in\D,
$$
where the minimum is taken over all curves $C$ in $\D$ from $z_1$ to $z_2$; the curve that minimizes the length is called the hyperbolic geodesic joining $z_1$ and $z_2$ (see \cite{Pom}, p.6). It is well known that, for $z_1, z_2\in\D$, the hyperbolic geodesic joining $z_1$ and $z_2$ is the arc of the circle through $z_1$ and $z_2$, orthogonal to $\partial {\D}$. The generalized Gehring-Hayman theorem states that, if $f$ is a conformal mapping of $\D$ and $J\subset\D$ is any Jordan arc from $z_1$ to $z_2$, then 
$$
\int_{S}|f'(z)||dz|\leq A\int_{J}|f'(z)||dz|,
$$
where $A$ is an absolute constant and $S$ is the hyperbolic geodesic joining $z_1$ and $z_2$. We refer to \cite[Theorem 4.20]{Pom} and \cite{Pom1} for more information on the generalized Gehring-Hayman theorem. Numerous significant results and conjectures have been derived from the above-mentioned Gehring-Hayman inequality, among which, one is the famous Hall's conjecture.
We add here that, despite remaining unsolved for an extended period, H\"{a}st\"{o} and Ponnusamy (see \cite{PP}) recently proved it. In 1968, Jaenisch \cite{J} improved the second part of Theorem~A by showing $A\geq 4.56$ with numerical computations and reducing the range of $A$ to $4.56\leq A \leq 17.45$. However, the disc version of the Theorem~A is not true if we merely take $f$ to be analytic in $\D$ and continuous on $\partial \D^{+}$. For example, the mapping 
$$
f_0(z) = e^{i\left(\frac{1+z}{1-z}\right)},~~z\in\D,
$$
is analytic in $\D$ and continuous on ${\partial\D}^{+}$. Here, $f_0$ maps ${\partial\D}^{+}$ onto a curve of length 1, but the length of $f_0((-1,1))$ is infinite. Also, Theorem~A cannot be extended to quasi-conformal mappings (see \cite[Theorem~4]{GH}). Thus, a natural question arises: Does Theorem A hold in the case of meromorphic univalent functions?

In this article, we give an answer to this question by proving a result similar to Theorem~A for meromorphic univalent functions with a simple pole in $\D$. For this purpose, we have chosen our $f$ to be a meromorphic and univalent function in $\D$ with a simple pole at $z=p\; (0<p<1)$ that remains continuous on $C'$. Let $\ell=\|f(C')\|$, $\ell'=\|f(I)\|$, where $\|\gamma\|$ denotes the length of $\gamma$. In Theorem~\ref{Th1}, we will prove the existence of a constant $A_p$ satisfying $\ell'\leq A_p\ell$ and provide a range for $A_p$. For the lower bound of $A_p$, consider the function 

$$
k_p(z)=\frac{pz}{(1-pz)(p-z)},~~z\in \D,
$$
which is meromorphic and univalent in $ \mathbb{D}$ with a simple pole at $z=p\,(0<p<1)$. It can be shown that (see \cite{AW})  
$$
k_p(\mathbb{D})=\hat{\mathbb{C}}\setminus\left[\frac{-p}{(1-p)^2},\frac{-p}{(1+p)^2}\right].
$$
Also, $k_p$ is continuous on $\partial \D$ and 
$$
k_p(-1)=\frac{-p}{(1+p)^2};~~ k_p(\pm i)=\frac{-p}{1+p^2}.
$$
Since $k_p$ maps $\partial \D$ onto the interval $[-{p}{(1-p)^{-2}},-{p}{(1+p)^{-2}}]$, we conclude that $k_p$ maps $C'$ onto the interval $(-{p}(1+p^2)^{-1},-{p}{(1+p)^{-2}}]$ twice. Thus,
\beq \label{p1eq1.1}
\|k_p(C')\|=2\bigg(\frac{p}{1+p^2}-\frac{p}{(1+p)^2}\bigg)=\frac{4p^2}{(1+p^2)(1+p)^2}.
\eeq 
Also, for $r \in (-1,1)$,
$$
k_p(ir)=\frac{ir}{(1-irp)(1-\frac{ir}{p})}
=\frac{-r^2(p+\frac{1}{p})+ir(1-r^2)}{(1-r^2)^2+r^2(p+\frac{1}{p})^2}.
$$
So, $k_p$ maps $I$ conformally onto the circle centered at $(-{p}/{2(1+{p^2})},0)$ with radius ${p}/{2(1+{p^2})}$  excluding the point $(-{p}(1+p^2)^{-1},0)$. Consequently,
\beq \label{p1eq1.2}
\|k_p(I)\|=\frac{2\pi}{2(p+\frac{1}{p})}=\frac{p\pi}{1+p^2}.
\eeq
Combining equations (\ref{p1eq1.1}) and (\ref{p1eq1.2}), we obtain
\be\nonumber
\frac{\|k_p(I)\|}{\|k_p(C')\|}=\frac{(1+p)^2 \pi}{4p}.
\ee
Therefore, we deduce that 
\be\label{p1eq1.3}
A_p\geq \frac{(1+p)^2\pi}{4p},~ p\in (0,1).
\ee

Following this, several lemmas will be demonstrated, which will be used to establish the existence of a constant $A_p$ in Theorem~$\ref{Th1}$. Furthermore, we will identify a specific region in $\D$ for the pole where the conclusion of Theorem~\ref{Th1} remains valid. Finally, we will give an extension of Theorem~$\ref{Th1}$ by replacing  vertical diameter $I$ with any hyperbolic geodesic in $\D$ that is symmetric about the real line and $C'$ with the corresponding arc of $\partial \D$ that passes through $-1$ and joins the endpoints of the aforementioned hyperbolic geodesic. This is the content of Theorem~$\ref{Th2}$.
\section{Preliminary Results}
Since we aim to determine 
$$
\ell'=\|f(I)\|=\int_{I}|f'(z)||dz|,
$$
we need an upper bound for the distortion $|f'(z)|$ at the points $z$ of $I$. Let $0<\alpha<1$ and $\gamma$ be the hyperbolic geodesic in $\D$ which is symmetric about the real axis and passes through the point ${\alpha}/{(1+\sqrt{1-{\alpha}^2})}$. Since hyperbolic geodesics of $\D$ are arcs of circles in $\D$ that are perpendicular to $\partial \D$, we obtain 
\be\label{Eq0.5}
\gamma = \{(x,y): x^2+y^2-({2}/{\alpha})x+1=0, x^2+y^2<1 \}.
\ee 
Now, $\gamma$ meets $\partial \D$ at the points $\alpha \pm i\sqrt{1-{\alpha}^2}$ and determines two hyperbolic half-planes; let $\Omega$ be the one that contains the origin. Then, the following function (see \cite[p. 85]{MM}) 
\beq \label{Eq1}
\phi_{\alpha}(z)=\frac{z(1-\alpha z)}{z-\alpha},\;z\in\Omega,
\eeq
maps $\Omega$ conformally onto $\D$ with $\phi_{\alpha}(0)=0$ and $\phi_{\alpha}'(0)=-{1}/{\alpha}$.
\blem \label{L1}
If $f$ maps $\Omega$ conformally onto a domain $D$, then 
\beq \nonumber
|f'(z)|\leq \frac{4|{\phi_{\alpha}}'(z)|}{1-|{\phi_{\alpha}(z)}|^2}\rho(f(z)),~~~~  z\in \Omega,
\eeq
where $\phi_{\alpha}$ is the function in $(\ref{Eq1})$ that maps $\Omega$ conformally onto $\D$ and $\rho(f(z))$ denotes the distance from $f(z)$ to the boundary of $D$.     
\elem
\bpf
We note that $f\circ {\phi_{\alpha}}^{-1}$ maps $\D$ conformally onto $D$. Applying the Koebe distortion theorem (see \cite[ p.9]{Pom}) to $f\circ {\phi_{\alpha}}^{-1} $, we get
\beq \label{L1Eq1}
|(f\circ {\phi_{\alpha}}^{-1})'(\zeta) |\leq\frac{4 \rho(f\circ {\phi_{\alpha}}^{-1}(\zeta))}{1-|\zeta|^2},\; \zeta\in\D.
\eeq
Now $|(f\circ {\phi_{\alpha}}^{-1})'(\zeta) |=|f'({\phi_{\alpha}}^{-1}(\zeta))|\;|({\phi_{\alpha}}^{-1})'(\zeta)|$.
Replacing $\zeta$ by $\phi_{\alpha}(z)$, we have 
$$
|(f\circ {\phi_{\alpha}}^{-1})'(\zeta) |=\frac{|f'(z)|}{|\phi_{\alpha}'(z)|}
$$
and 
$$
\rho(f\circ {\phi_{\alpha}}^{-1}(\zeta))=\dist({f\circ {\phi_{\alpha}}^{-1}}(\zeta), \partial({f\circ {\phi_{\alpha}}^{-1}}(\D)))=\dist(f(z), \partial f(\Omega)) =\rho(f(z)).
$$
Thus, letting $\zeta=\phi_{\alpha}(z)$ in (\ref{L1Eq1}), we get the desired inequality.
\epf

The above lemma provides a local measure of distortion of a conformal mapping defined on $\Omega$. Now, we will proceed toward our main goal, i.e. proving Theorem~\ref{Th1} using a similar method of proof which was adopted by Gehring and Hayman in \cite{GH}. Let
\beq \label{Eq5}
g(z)=\frac{-1-iz}{z+i},~z\in\D.
\eeq
It is easy to see that, $g$ maps $\D$ conformally onto $\mathbb{H}$ with 
\be\label{Eq5.25}
g(0)=i ~~~\text{ and }~~~g\left(\frac{\alpha}{1+\sqrt{1-{\alpha}^2}}\right)=-\alpha+i\sqrt{1-\alpha^2}.
\ee
Also, $g(I)=\mathbb{I}^{+}$ and $g(C')={\IR}^{+}$. Since $g$ is a M\"{o}bius map, it will map geodesics onto geodesics, and geodesics of $\mathbb{H}$ are semi-circles orthogonal to the real axis and vertical half-lines. Now, $g$ maps $\gamma$ which is defined in ($\ref{Eq0.5}$) conformally onto the semi-circle 
\be\label{Eq5.5}
\sigma = \{(x,y): x^2+y^2+({2}/{\alpha})x+1=0, y>0\}.
\ee
It is clear that $\sigma$ is perpendicular to the real line. Thus, $\sigma$ is a hyperbolic geodesic in $\mathbb{H}$. Let $\Omega_1$ be the hyperbolic half-plane of $\mathbb{H}$ determined by $\sigma$ that contains the point $i$. Therefore, $g$ maps $\Omega$, which is defined in Lemma~$\ref{L1}$, conformally onto $\Omega_1$. The following lemma is the upper half-plane version of Lemma~$\ref{L1}$.
\blem \label{L2}
Let $w$ map $\Omega_1$ conformally onto a domain $D$. Then 
$$
|w'(z)|\leq \frac{8|{\phi_{\alpha}}'(-\frac{iz+1}{z+i})|}{|z+i|^2(1-|{\phi_{\alpha}}(-\frac{iz+1}{z+i})|^2)}\rho(w(z)),~~~~~z\in\Omega_1,
$$
where $\rho(w(z))$ denotes the distance from $w(z)$ to the boundary of $D$ and $\phi_\alpha$ is as in $(\ref{Eq1})$.
\elem
\bpf
Since $w\circ g$ maps $\Omega$ conformally onto $D$, from Lemma $\ref{L1}$ we get
\be \label{L2Eq1}
|(w\circ g)'(\zeta)|\leq \frac{4|{\phi_{\alpha}}'(\zeta)|}{1-|\phi_{\alpha}(\zeta)|^2}\rho({w\circ g}(\zeta)).
\ee
Now, letting $\zeta=g^{-1}(z)=-{(iz+1)}/{(z+i)}$, we have
$$
\rho({w\circ g(\zeta)})=\dist(w\circ g(\zeta), \partial(w\circ g(\Omega)))
=\dist(w(z), \partial w(\Omega_1))
=\rho(w(z))
$$
and
$$
\frac{2|w'(z)|}{|\frac{-iz-1}{z+i}+i|^2}\leq \frac{4|{\phi_{\alpha}}'(\frac{-iz-1}{z+i})|}{1-|\phi_{\alpha}(\frac{-iz-1}{z+i})|^2}\rho (w(z))\\\nonumber
$$
from (\ref{L2Eq1}). The desired inequality follows from the above estimate.
\epf

The lemma below provides an upper bound for $\rho (w)$ (see \cite[Lemma~2]{GH}).

\vspace{0.2cm}
\noindent\textbf{Lemma~A.}\label{L3}
\textit{Suppose that $M$ is a simply connected domain and $\beta$ is a set of boundary points of $M$ which lie within distance $d$ of a point $w_0$ not in $M$. Then if $u(w)$ is harmonic in $M$, satisfies $0\leq u(w)\leq 1$ there, and approaches $0$ as $w$ approaches any boundary points of $M$ not in $\beta$, we have} 
$$
\rho(w)\leq d\,\cot^2\left(\frac{\pi}{4}u(w)\right).
$$

We now introduce the notion of harmonic measure, which will be used to prove Lemma~$\ref{L4}$ below. Given a domain $D$, a point $z \in D$, and a Borel set $E$ of $\overline{D}$, where $\overline{D}$ denotes the closure of $D$, $\omega(z,E,D)$ denotes the harmonic measure at $z$ of $\overline{E}$ with respect to the component of $D\setminus \overline{E}$ containing $z$. The function $\omega(.,E,D)$  precisely represents the solution to the generalized Dirichlet problem with boundary data $\phi=1_{E}$ (see \cite[ch.3]{LV}, \cite[p.297]{CC}). Clearly, $\omega(z,E,D)$ is a harmonic function satisfying $0\leq\omega(z,E,D)\leq1$. The next theorem provides a relation between the harmonic measures defined on a domain and on its subdomain (see \cite[p.282]{Bet} and \cite[p.302]{CC}).

\vspace{0.2cm}
\noindent\textbf{Theorem B.}\label{Th2}
\textit{Let $D_1$ and $D_2$ be two domains in $\IC$. Assume that $D_1\subset D_2$ and let $F\subset \partial D_2$ be a closed set. Let $\eta=\partial D_1 \setminus \partial D_2$. Then for $z\in D_1$},
\beq \nonumber
\omega(z,F,D_2)=\omega(z,F,D_1)+\int_{\eta}\omega(z,ds,D_1)\omega(s,F,D_2).
\eeq
After all these preparations, we prove the following lemma which will be used to establish Theorem~$\ref{Th1}$.
\blem \label{L4}
Let $\alpha\in (1/\sqrt{2},1)$ and $\Omega_1$ be the same as in Lemma~$\ref{L2}$. Suppose that $0<a<b$ and that $\beta$ and $\beta'$ are the segments joining $a$ to $b$ and $ia$ to $ib$, respectively. If $w$ is a conformal mapping of $\Omega_1$, which remains continuous on $\beta$ and maps $\beta$ and $\beta'$ onto curves of length $\ell$ and $\ell'$, then we have 
$$
\ell'< B_\alpha\ell,
$$
where $B_\alpha$ is a constant which depends on $\alpha$ and ${b}/{a}$.
\elem
\bpf
Let $\omega(z,\beta,\mathbb{H})$ and $\omega(z,\beta,\Omega_1)$ denote the harmonic measures at the point $z$ of $\beta$ with respect to $\mathbb{H}$ and $\Omega_1$, respectively. From Theorem~B, we get
\beq \label{L3Eq1}
\omega(z,\beta,\Omega_1)=\omega(z,\beta,\mathbb{H})-\int_{\sigma}\omega(z,ds,\Omega_1) \omega(s,\beta,\mathbb{H}),
\eeq
where $\sigma$ is given in ($\ref{Eq5.5}$).
Note that $\omega(z,\beta,\mathbb{H})$ is ${1}/{\pi}$ times the angle which $\beta$ subtends at $z$. Hence, for $z\in\beta'$ (see \cite[Lemma~3]{GH}),
\be \label{L3Eq2}
\omega(z,\beta,\mathbb{H}) \geq \frac{1}{\pi} \tan^{-1} \left( \frac{b-a}{b+a}\right).
\ee
Next, we will calculate an upper bound of $\omega(s,\beta,\mathbb{H})$ on $\sigma$. Let $s=(x,\sqrt{r^2-(x+m)^2})$ be an arbitrary point on $\sigma$, where 
$$
m=\frac{1}{\alpha},~~ r=\frac{\sqrt{1-\alpha^2}}{\alpha}~~\text{ and}~~ 
\frac{-1-\sqrt{1-\alpha^2}}{\alpha}\leq x\leq\frac{-\alpha}{1+\sqrt{1-\alpha^2}}.
$$ 
Then, 
\beq\nonumber
\omega(s,\beta,\mathbb{H})&=&\frac{1}{\pi}\left[\tan^{-1}\frac{b-x}{\sqrt{r^2-(x+m)^2}}-\tan^{-1}\frac{a-x}{\sqrt{r^2-(x+m)^2}}\right]\\\nonumber
&=&\frac{1}{\pi}\tan^{-1}\frac{(b-a)\sqrt{r^2-(x+m)^2}}{ab-1-2mx-(b+a)x}.
\eeq
Since $ab-1-2mx>0$ and $-(a+b)x>0$ for all $x\in[-({1+\sqrt{1-\alpha^2}})/{\alpha},-{\alpha}/({1+\sqrt{1-\alpha^2}})]$ and $\tan^{-1}(x)$ is an increasing function, we obtain from the equation above 
\beq\label{L3Eq3}
\omega(s,\beta,\mathbb{H})<\frac{1}{\pi}\tan^{-1}\left(\frac{b-a}{b+a}\psi(x)\right),
\eeq
where
$$
\psi(x)=\frac{-\sqrt{r^2-(x+m)^2}}{x}.
$$
Here, $\psi$ is continuous in $[-({1+\sqrt{1-\alpha^2}})/{\alpha},-{\alpha}/({1+\sqrt{1-\alpha^2}})]$ and  attains its maximum at the point $x=-\alpha$. Then from (\ref{L3Eq3}), we get
$$
\omega(s,\beta,\mathbb{H})<\frac{1}{\pi} \tan^{-1} \left( \frac{b-a}{b+a}\psi(-\alpha)\right)=\frac{1}{\pi}\tan^{-1}\left(\frac{b-a}{b+a}\frac{\sqrt{1-\alpha^2}}{\alpha}\right)
$$
for $s\in\sigma$.
Therefore,
\beq \nonumber
\int_{\sigma}\omega(z,ds,\Omega_1)\omega(s,\beta,\mathbb{H})
&<&\frac{1}{\pi} \tan^{-1} \left(\frac{b-a}{b+a} \frac{\sqrt{1-\alpha^2}}{\alpha}\right) \int_{\sigma} \omega(z,ds,\Omega_1)\\\nonumber
&=&\frac{1}{\pi} \tan^{-1} \left(\frac{b-a}{b+a} \frac{\sqrt{1-\alpha^2}}{\alpha}\right)\omega(z,\sigma,\Omega_1)\\ \label{L3Eq4}
&\leq&\frac{1}{\pi} \tan^{-1} \left(\frac{b-a}{b+a} \frac{\sqrt{1-\alpha^2}}{\alpha}\right).
\eeq
Thus, from (\ref{L3Eq1}), using (\ref{L3Eq2}) and (\ref{L3Eq4}), for $z\in\beta'$, we have 

\beq \label{L3Eq5}
&&\omega(z,\beta,\Omega_1)>\frac{\xi}{\pi},
\eeq
where
$$
\xi=\tan^{-1} \left( \frac{b-a}{b+a}\right)-\tan^{-1} \left(\frac{b-a}{b+a} \frac{\sqrt{1-\alpha^2}}{\alpha}\right).
$$
Since $\alpha \in (1/\sqrt{2},1)$, consequently $0<{\sqrt{1-\alpha^2}}/{\alpha}<1$. Thus, $0<\xi<\pi/4$.
Suppose $\beta_1=w(\beta)$ and $\beta'_1=w(\beta')$. Let $w_0$ denote the point on $\beta_1$ that divides its arclength in half. Define 
$$
u(w)=\omega(z(w),\beta,\Omega_1),
$$
where $z(w)$ is the inverse of $w$. Therefore, $\beta_1$ lies in the disk $|w-w_0|\leq d=\ell/2$ and $u(w)$ and $\beta_1$ satisfy the hypothesis of Lemma~A. Now, (\ref{L3Eq5}) implies
$$
u(w)>\frac{\xi}{\pi}
$$
at each point $w\in\beta'_1$  and we conclude from Lemma A and Lemma~\ref{L3} that 
$$
|w'(z)|< \frac{4\ell|{\phi_{\alpha}}'(-\frac{iz+1}{z+i})|}{|z+i|^2(1-|\phi_{\alpha}(-\frac{iz+1}{z+i})|^2)} \cot^2\left(\frac{\xi}{4}\right)
$$
at each point $z\in\beta'$. By integrating the above inequality along $\beta'$, we obtain
\beq\nonumber
\ell'=\int_{\beta'}|w'(z)||dz|&<&4\ell\;\cot^2\left(\frac{\xi}{4}\right)\int_{[ia,ib]}\frac{|{\phi_{\alpha}}'(-\frac{iz+1}{z+i})|}{|z+i|^2(1-|\phi_{\alpha}(-\frac{iz+1}{z+i})|^2)}|dz|\\\nonumber &=& 4\ell\; \cot^2\left(\frac{\xi}{4}\right)\int_{a}^{b}\frac{|{\phi_{\alpha}}'(i(\frac{1-y}{1+y}))|}{(1+y)^2(1-|{\phi_{\alpha}}(i(\frac{1-y}{1+y}))|^2)}dy.\\\nonumber
\eeq
By letting $r=(1-y)/(1+y)$ and using (\ref{Eq1}), from the above inequality we get
\beq \nonumber
\ell'&<& 2\ell\; \cot^2\left(\frac{\xi}{4}\right)\int_{\frac{1-b}{1+b}}^{\frac{1-a}{1+a}}\frac{|{\phi_{\alpha}'(ir)}|}{1-|\phi_{\alpha}(ir)|^2}dr\\\nonumber
&=&2\ell\; \cot^2\left(\frac{\xi}{4}\right)\int_{\frac{1-b}{1+b}}^{\frac{1-a}{1+a}}\frac{\sqrt{(1-r^2)^2+4\alpha^2 r^2}}{\alpha(1-r^4)}dr\\\nonumber
&<&2\ell\; \cot^2\left(\frac{\xi}{4}\right)\int_{\frac{1-b}{1+b}}^{\frac{1-a}{1+a}}\frac{\sqrt{(1-r^2)^2+4r^2}}{\alpha(1-r^4)}dr\\\nonumber
&=&2\ell\; \cot^2\left(\frac{\xi}{4}\right)\int_{\frac{1-b}{1+b}}^{\frac{1-a}{1+a}}\frac{dr}{\alpha(1-r^2)}.
\eeq
Thus, we have
$$
\ell'<B_\alpha\;\ell,
$$
where 
\be \nonumber
B_\alpha=\frac{1}{\alpha} \cot^2\left(\frac{\xi}{4}\right)\log\left(\frac{b}{a}\right).
\ee
Note that $\xi$ depends on $\alpha$ and ${b}/{a}$. Therefore, for $\alpha\in({1}/{\sqrt{2}},1)$, $B_\alpha$ is a constant that depends on $\alpha$ and ${b}/{a}$.
\epf
\section{Main Results}
We now state and prove the main result of this article.
\bthm\label{Th1}
Let $p\in (\sqrt{2}-1,1)$. Suppose $f$ is a meromorphic univalent function in $\D$ with a simple pole at $z=p$ that remains continuous on the semi-circle $C'$. If $\ell=\|f(C')\|$ and $\ell'=\|f(I)\|$, then there exists a constant $A_p$ such that
\be\label{Th1Eq1}
\ell'\leq A_p\;\ell
\ee
with 
\be\label{Th1Eq2}
\frac{(1+p)^2 \pi}{4p}\leq A_p\leq \displaystyle{\min_{q\in(1,\infty)}M_p (q)},
\ee
where
\be\label{Th1Eq3}
M_p(q)=\frac{(1+p^2)}{2p}\cot^2\left(\frac{\tan^{-1}\left(\frac{q-1}{q+1}\right)-\tan^{-1}\left(\frac{(1-p^2)(q-1)}{2p(q+1)}\right)}{4}\right)\log q.
\ee
\ethm

\bpf
Let $g$ be the mapping defined in $(\ref{Eq5})$. Let
\be\label{Th1Eq3.5}
p=\frac{\alpha}{1+\sqrt{1-{\alpha}^2}},\;\text{ i.e. }\;\alpha=\frac{2p}{(1+p^2)}.
\ee 
It is easy to see that $p\in(\sqrt{2}-1,1)$ if and only if $\alpha\in({1}/{\sqrt{2}},1)$. Then $w=f\circ g^{-1}$ is meromorphic and univalent in $\mathbb{H}$ with a simple pole at $\alpha':= g(p)=-\alpha+i\sqrt{1-\alpha^2}$ (see (\ref{Eq5.25})) that remains continuous on ${\IR}^{+}$. Let $w_1=w|_{\Omega_1}$, where $\Omega_1$ is the same as in Lemma~$\ref{L3}$. Since $\alpha'$ lies on $\sigma$, which is defined in ($\ref{Eq5.5}$), then $w_1$ is conformal on $\Omega_1$ and remains continuous on ${\IR}^{+}$. Choose $1<q<\infty$. For $n=0,\pm1,\pm2,...$, let $\beta_n$ and $\beta'_n$ denote the segments joining $q^n$ to $q^{n+1}$ and $iq^n$ to $iq^{n+1}$, respectively. Let $\ell_n$ and $\ell'_n$ denote the lengths of these segments under $w$ (or $w_1$). Then from Lemma~\ref{L4}, we get 
\be\label{Th1eq3.75}
\ell'_n<B_\alpha(q)\ell_n,
\ee
where 
\be \label{Th1Eq4}
B_\alpha(q)=\frac{1}{\alpha} \cot^2\left(\frac{\tan^{-1}\left(\frac{q-1}{q+1}\right)-\tan^{-1}\left(\frac{\sqrt{1-\alpha^2}(q-1)}{\alpha(q+1)}\right)}{4}\right)\log q.
\ee
Now $$\ell=\|f(C')\|=\|(f\circ g^{-1})({\IR}^{+})\|=\|w_1({\IR}^{+})\|
$$
and 
$$
\ell'=\|f(I)\|=\|(f\circ g^{-1})(\mathbb{I}^{+})\|=\|w_1(\mathbb{I}^{+})\|.
$$
Thus, from (\ref{Th1eq3.75}), we get
\be \label{Th1Eq5}
\ell'=\sum_{n=-\infty}^{\infty}\ell'_n\leq B_\alpha(q)\sum_{n=-\infty}^{\infty}\ell_n=B_\alpha(q)\ell.
\ee
 Now in (\ref{Th1Eq4}) and (\ref{Th1Eq5}), letting $\alpha=  {2p}/{(1+p^2)}$, we obtain
$$
\ell'\leq M_p (q)\ell, 
$$
where $M_p (q)=B_{{2p}/{(1+p^2)}}(q)$, stated in ($\ref{Th1Eq3}$),
is a constant depending on $q\in(1,\infty)$. Clearly, $M_p (q)>0$ and continuous in the interval $(1,\infty)$. Moreover, $M_p (q)\to +\infty$ as $q\to 1^{+} \text{ or } q\to+\infty$. So, $M_p (q)$ attains its minimum in $(1,\infty)$. Thus, from the previous inequality, we obtain
$$
\ell'\leq\left(\displaystyle{\min_{q\in(1,\infty)}M_p (q)}\right) \ell.
$$
This proves the existence of a constant $A_p$ satisfying the inequality ($\ref{Th1Eq1}$) and the right-hand side inequality of ($\ref{Th1Eq2}$). The left-hand side inequality of ($\ref{Th1Eq2}$) has already been established in (\ref{p1eq1.3}). Hence, the theorem follows.
\epf

In Table~\ref{table1}, we demonstrate ranges of $A_p$ for some specific values for $p \in (\sqrt{2}-1,1)$ with the help of Mathematica 11.0. Additionally, we indicate the point $q_p$, where the minimum value of $M_p (q)$ is achieved.
\begin{table}[ht]
	
	\begin{tabular}{|c |c
			|c |c |}		
		\hline 
		$p$ & $ q_p$ &  $\displaystyle{\min_{q\in(1,\infty)}M_p}$ &  Range of $A_p$   \\ [0.1ex] 
		\hline
		0.999 & 5.55 & 73.42 & (3.14, 74)\\
		\hline
		0.99  &  5.52  & 74.99 & (3.14,  75)\\
		\hline
		0.9  &  5.19  & 95.49 & (3.15, 96) \\
		\hline
		0.8 & 4.78   & 135.73 & (3.18, 136) \\
		\hline
		0.7 &  4.33  & 221.80 &  (3.24, 222)   \\
		\hline
		0.6 &  3.85  & 471.01 & (3.35, 472) \\
		\hline
		0.5 &  3.37  & 1984.43 & (3.53, 1985) \\
		\hline
		0.45 & 3.13 & 10811.10 & (3.66, 10812) \\
		\hline
		0.423 & 3.01 & 174258 & (3.75, 174259) \\
		\hline
	\end{tabular} 
	\vspace{0.2cm}\caption{Ranges of $A_p$ for various values of $p\in$ ($\sqrt{2}-1,1$)}\label{table1}
\end{table}

\brs 
From Theorem~\ref{Th1}, we have the following immediate consequences.
\begin{itemize}
	\item [(i)] Letting $p\to1^{-}$ in Theorem~\ref{Th1}, we have
	$$
	\ell'\leq A_{1}\ell \;\mbox{ with }\;\pi\leq A_{1}<74,
	$$
	which coincides with the case of conformal mappings given in Theorem~A.
	\item [(ii)] In Theorem~$\ref{Th1}$, if the simple pole of $f$ lies anywhere on $\gamma$,  defined in ($\ref{Eq0.5}$), not necessarily at $z=p={\alpha}/({1+\sqrt{1-{\alpha}^2}})$, then it is easy to see from the proof of Theorem~$\ref{Th1}$ that one will get the same  constant $A_p$ for which (\ref{Th1Eq1}) and ($\ref{Th1Eq2}$) hold.
	\item [(iii)] Let $\gamma_0$ be the hyperbolic geodesic in $\D$ that is symmetric about the real axis and passes through the point $\sqrt{2}-1$. Using ($\ref{Th1Eq3.5}$), we see from ($\ref{Eq0.5}$) that
	\be\label{R3Eq1}
	\gamma_0 = \{(x,y): x^2+y^2-2\sqrt{2}x+1=0, x^2+y^2<1 \}.
	\ee
Denote $\Omega_0$ the hyperbolic half-plane determined by $\gamma_0$ that does not contain the origin. Suppose $h$ is a meromorphic univalent function in $\D$ with a simple pole at $\Tilde{p}$ in $\Omega_0$ that is also continuous on $C'$. Let us assume that the hyperbolic geodesic passing through $\Tilde{p}$ and symmetric about the real axis meets the real line at the point $\Tilde{p_1}$. Then, $\sqrt{2}-1<\Tilde{p_1}<1$. Therefore, from Remark (ii), we have
	$$ 
	{\ell_1}'\leq A_{\Tilde{p_1}}\ell_1,
	$$ 
	where $\ell_1=\|h(C')\|$, ${\ell_1}'=\|h(I)\|$, and $A_{\Tilde{p_1}}$ is a constant satisfying ($\ref{Th1Eq2}$), where we replace $p$ with ${\Tilde{p_1}}$. 
	
\end{itemize}
\ers

We believe that for the rest of the  values of $p\in (0,1)$, i.e. for $p\in (0, \sqrt 2 -1]$, the conclusion of Theorem~1 is true. Therefore, we  conjecture the following:
\bcon
Let $p\in (0, \sqrt 2 -1]$. Suppose $f$ is a meromorphic univalent function in $\D$ with a simple pole at $p$ and satisfies the hypothesis of Theorem~1. Then the estimate $(\ref{Th1Eq1})$ is valid for some constant $A_p$.
\econ

Next, setting the pole $p$ in the interval $(\sqrt{2}-1,1)$, we wish to see the extent to which we can rotate the vertical diameter $I$ and the corresponding semi-circle $C'$ while maintaining the validity of the conclusion stated in Theorem~$\ref{Th1}$. This is the content of the following corollary.
\bcor 
Let $p\in(\sqrt{2}-1,1)$ and $|{\theta}|<\tan^{-1}({\sqrt{6p^2-p^4-1}}/{(1+p^2)})$. Let $I_1$ be a diameter of $\D$ which makes an angle $\theta$ with $I$ and $C_1$ be the arc of $\partial\D$ passing through the point $-1$ and joining the endpoints of $I_1$. Suppose $f$ is a meromorphic univalent function in $\D$ with a simple pole at $z=p$ and continuous on $C_1$. 
If $\ell=\|f(C_1)\|$ and $\ell'=\|f(I_1)\|$, then there exists a constant $A_{p,\theta}$ (depending on $p$ and $\theta$) such that $\ell'\leq A_{p,\theta}\,\ell$.
\ecor
\bpf
If $p\in(\sqrt{2}-1,1)$, then $|\theta|<\pi/4$.
Without loss of generality, let $\theta$ be positive. Now, rotating the unit disk by the angle $\theta$ in the clockwise direction, we see that  $I_1$ and $C_1$ are mapped onto $I$ and $C'$, respectively, and the pole $p$ has been moved to the point $pe^{-i\theta}$. A straightforward computation yields 
\be\label{p1eq5a}
\{z: |z|=p\}\cap \gamma_0=\{{(1+p^2)}/{2\sqrt{2}},\pm {\sqrt{6p^2-p^4-1}}/{2\sqrt{2}}\},
\ee
where $\gamma_0$ is the hyperbolic geodesic in $\D$ that is symmetric about the real axis and passes through the point $\sqrt{2}-1$, defined in  ($\ref{R3Eq1}$). By the hypothesis 
$$
\theta<\tan^{-1}({\sqrt{6p^2-p^4-1}}/{(1+p^2)}),
$$
i.e. the angle $\theta$ is less than the modulus of the arguments of the points indicated in (\ref{p1eq5a}), therefore we see that $pe^{-i\theta}\in\Omega_0$. Here,  $\Omega_0$ is the hyperbolic half-plane determined by $\gamma_0$ that does not contain the origin. Now, $f(ze^{i\theta})$ is meromorphic univalent in $\D$ with a simple pole at $pe^{i\theta}$ in $\Omega_0$ and remains continuous on $C'$. Hence, applying the same reasoning that was outlined in Remark~(iii)  to the function $f(ze^{i\theta})$, we see that, a constant $A_{p,\theta}$ exists (depending on $p$ and $\theta$) such that $\ell'\leq A_{p,\theta}\,\ell$ holds.
\epf

In Theorem~$\ref{Th1}$, we have specifically considered the vertical diameter $I$ and the associated arc of $\partial \D$ passing through the point $-1$. Now, we aim to extend Theorem~$\ref{Th1}$ to any hyperbolic geodesic symmetric about the real line, and its corresponding arc of $\partial \D$ passing through the point $-1$. To accomplish this, we select $\alpha_1$ on $\partial \D^{+}$ and designate $\gamma_1$ as the hyperbolic geodesic that is symmetric about the real axis and passes 
through the point $\alpha_1$. Let $\gamma_2$ be the arc of $\partial\D$ which passes through the point $-1$ and connects the endpoints of $\gamma_1$. Let
\be\label{Eq6}
p_0=\frac{1+\sqrt{2}\,{\rm{Re}}\,\alpha_1}{\sqrt{2}+{\rm{Re}}\,((1-i)\alpha_1)}.
\ee
It is easy to see that $p_0\in (-1,1)$.
In the following, we present a generalization of Theorem~$\ref{Th1}$.
\bthm
Let $p_1\in (p_0,1)$, where $p_0$ is given in $(\ref{Eq6})$. Suppose $f$ is a meromorphic univalent function in $\D$ with a simple pole at $p_1$ that remains continuous on $\gamma_2$. If $f$ maps $\gamma_2$ and $\gamma_1$ onto curves of lengths $\ell$ and $\ell'$, respectively, then 
$$
\ell'\leq A_{p'_1}\;\ell,
$$
where $A_{p'_1}$ is a constant satisfying the inequality $(\ref{Th1Eq2})$ with the substitution of $p$ with $p'_1$, where 
$$
p'_1=\frac{1-p_1\alpha_1+i (\alpha_1-p_1)}{p_1-\alpha_1+i(p_1\alpha_1 -1)}.
$$
\ethm
\bpf
Let
$$
g_1(z)=\frac{z-\alpha_1}{\alpha_1 z-1},~~~z\in\D.
$$
Then $g_1$ maps $\D$ conformally onto $\mathbb{H}$ in such a way that $\gamma_1$ and $\gamma_2$ are mapped onto $\mathbb{I}^{+}$ and ${\IR}^{+}$, respectively. Considering the map $g$ defined in ($\ref{Eq5}$) we see that the function $g^{-1}\circ g_1$ is a conformal mapping of $\D$ onto itself and maps $\gamma_1$ and $\gamma_2$ onto $I$ and $C'$, respectively. Also, for $r\in(-1,1)$, $|g_1(r)|=1$ and $|g(r)|=1$. Thus, $g^{-1}\circ g_1$ maps the interval $(-1,1)$ conformally onto itself. Since
\be\nonumber
(g^{-1}\circ g_1)^{-1}(\sqrt{2}-1)={g_1}^{-1}\left(\frac{-1+i}{\sqrt{2}}\right)
=p_0,
\ee
where $p_0$ is given in ($\ref{Eq6}$)
and $(g^{-1}\circ g_1)(1)=1$, therefore the interval $(p_0,1)$ is conformally mapped onto the interval $(\sqrt{2}-1,1)$ by $g^{-1}\circ g_1$. Thus, $f\circ (g^{-1}\circ g_1)^{-1}$ is a meromorphic univalent function in $\D$ with a simple pole at $(g^{-1}\circ g_1)(p_1)\in(\sqrt{2}-1,1)$ that remains continuous on $C'$ and maps $I$ and $C'$ onto curves of length $\ell'$ and $\ell$, respectively. Denote 
$$
p'_1=(g^{-1}\circ g_1)(p_1)=\frac{1-p_1\alpha_1+i (\alpha_1-p_1)}{p_1-\alpha_1+i(p_1\alpha_1 -1)}.
$$
Hence, the theorem follows by an application of Theorem~$\ref{Th1}$ to the function $f\circ (g^{-1}\circ g_1)^{-1}$, where we replace $p$ with $p'_1$.
\epf

\section{Acknowledgements}

\vspace{.4cm}

\noindent{\bf Statements and Declarations}: Nil.\\
\noindent{\bf Competing interests}: We have nothing to declare.\\
\noindent{\bf Conflict of interest}: Nil.\\

\end{document}